\documentclass[12pt]{amsart}

\usepackage{amsmath,amsthm,amssymb}
\font\teneufm=eufm10
\font\seveneufm=eufm7
\font\fiveeufm=eufm5
\newfam\frakturfam

\textfont\frakturfam=\teneufm
\scriptfont\frakturfam=\seveneufm
\scriptscriptfont\frakturfam=\fiveeufm

\newtheorem{proposition}{Proposition}

\newtheorem{lemma}{Lemma}
\newtheorem{theorem}{Theorem}
\newtheorem{corollary}{Corollary}

\newtheorem{problem}{Problem}
\def\bee{\begin{eqnarray}}
\def\bes{\begin{eqnarray*}}
\def\eee{\end{eqnarray}}
\def\ees{\end{eqnarray*}}

\def\a{\alpha}

\def\g{\gamma}

\def\t{\tau}
\def\d{\partial}
\def\Proof{{\sl Proof.}\ }

\title{Automorphisms of free metabelian Lie algebras, II}

\begin{document}
\date{}

\maketitle

\begin{center}

{\bf Ualbai Umirbaev}\footnote{Department of Mathematics,
 Wayne State University,
Detroit, MI 48202, USA; Center for Theoretical and Applied Algebra, 
Al-Farabi Kazakh National University, Almaty, 050040, Kazakhstan; 
and Institute of Mathematics and Mathematical Modeling, Almaty, 050010, Kazakhstan,
e-mail: {\em umirbaev@wayne.edu}}

\bigskip

\dedicatory{In memory of Professor V.A. Roman'kov (1948--2023)}

\end{center}

\begin{abstract} 
We show that every automorphism of the free metabelian Lie algebra $M_n$ of rank $n\geq 4$ over an arbitrary field $K$ is almost tame; that is, it can be expressed as a product of Chein automorphisms (also called one-row transformations). Moreover, we prove that, for $n\geq 4$, the group  $\mathrm{Aut}(M_n)$ of all automorphisms of the free metabelian Lie algebra $M_n$ over a field $K$ of characteristic not equal to $3$ is generated by its linear automorphisms, together with a single quadratic and a single cubic automorphism. The same result holds for fields of any characteristic if $n\geq 5$.  
\end{abstract}



\tableofcontents

\section{Introduction}

\hspace*{\parindent}

There are many interesting and profound results concerning the automorphism group of the free metabelian Lie algebra $M_n$ of rank $n\geq 2$ over a field $K$. Some of them are described in \cite{AutoI}. In this paper, we continue the study of tame and almost tame automorphisms of the free metabelian Lie algebra $M_n$ of rank $n\geq 4$ over $K$, initiated in \cite{AutoI}. Furthermore, we use the definitions and results of \cite{AutoI} throughout and refer the reader there for details.

Free metabelian Lie algebras admit non-elementary automorphisms that change only one variable and fix the rest. Chein proved that they give examples of wild automorphisms for free metabelian groups of rank $3$ \cite{Chein}. All automorphisms of $M_n$ of this type (including elementary ones) are called {\em Chein automorphisms} or {\em one-row transformations} \cite{GGR,Romankov95,AutoI}. It was shown in \cite{AutoI} that all Chein automorphisms of lower degree at least 4, as well as all exponential automorphisms of lower degree at least 5, of the free metabelian Lie algebra $M_n$ of rank $n\geq 4$ over a field of characteristic different from 3, are tame.  The same results hold for fields of any characteristic when $n\geq 5$. Moreover, it was shown that a large class of automorphisms (see \cite[Theorem 3]{AutoI} for the precise formulations) of $M_n$ of rank $n\geq 4$ that move only two variables are {\em almost tame} \cite{Nauryzbaev09}; that is, they can be expressed as a product of Chein automorphisms.   

In 2015, Nauryzbaev \cite{Nauryzbaev15} proved that the group of tame automorphisms of $M_n$ is generated by all linear automorphisms and the quadratic automorphism 
\bee\label{f1}
(x_1+[x_2,x_3],x_2,\ldots,x_n)
\eee
if either $n\geq 4$ and $K$ is a field of characteristic $\neq 3$, or $n\geq 5$ and $K$ is a field of arbitrary characteristic.

 The first example of a non-elementary Chein automorphism is 
\bee\label{f2}
(x_1+[[x_2,x_3],x_1],x_2,\ldots,x_n). 
\eee
  This automorphism is wild if $n=3$ \cite{Nauryzbaev09}, and is the first candidate for wildness for all $n\geq 4$. 
 
 \begin{problem}\label{pr1}
 	Is the Chein automorphism (\ref{f2}) of the free metabelian algebra Lie algebra $M_n$ of rank $n\geq 4$ tame? 
 \end{problem}

Recall that the well known Bachmuth-Mochizuki-Roman'kov Theorem \cite{BM85,Romankov85} states that every automorphism of the free metabelian group of rank $n\geq 4$ is tame. 

In this paper we show that every automorphism of $M_n$ of rank $n\geq 4$ is almost tame. 
This is an exact analogue of the above-mentioned Bachmuth–Mochizuki–Roman'kov Theorem modulo Problem \ref{pr1}. More precisely, if Problem \ref{pr1} has a positive solution, then an exact analogue of the Bachmuth–Mochizuki–Roman'kov Theorem holds for free metabelian Lie algebras.

We also show that, if $K$ is a field of characteristic different from $3$, then the automorphism group $\mathrm{Aut}(M_n)$ of the free metabelian Lie algebra $M_n$ of rank $n\geq 4$  is generated by all linear automorphisms together with the quadratic automorphism (\ref{f1}) and the cubic automorphism (\ref{f2}). If $n\geq 5$ then the same result holds for fields of any characteristic.

Thus the question about the existence of wild automorphisms for the free metabelian Lie algebra $M_n$ of rank $n\geq 4$ is still open as well as in the case polynomial and free associative algebras. But in this case it is reduced to the concrete automorphism (\ref{f2}).

The paper is organized as follows. In Section 2, we recall some notation and definitions from \cite{AutoI}. In Section 3,  we establish matrix notations for automorphisms using the Jacobian matrices, and give matrix formulations of some results of \cite{AutoI}.  
Section 4 is devoted to the study of some special class of matrices that includes the Jacobian matrices of an important class of automorphisms. In Section 5, we prove a lemma on the elimination of a variable under certain conditions.  Section 6 is devoted to the study of one exceptional automorphism 
\bee\label{f3}
(x_1+[[x_1,x_2],x_1],x_2+[[x_1,x_2],x_2],x_3,\ldots,x_n)
\eee
of $M_n$.  In Section 7, we show how to eliminate a variable in the general case. 
In Section 8, the main result and some of its corollaries are proved.

\section{Preliminaries}

\hspace*{\parindent}

Let $M_n$ be the free metabelian Lie algebra over a field $K$ with free generators $x_1,\ldots,x_n$. 
The set of all right normed Lie monomials of the form 
\bee\label{f4}
[\ldots[[x_i,x_{i_1}],x_{i_2}],\ldots,x_{i_k}], \ i>i_1\leq i_2\leq\ldots \leq i_k, 
\eee
is a linear basis of $M_n'=[M_n,M_n]$ (see, for example \cite[p. 139]{Bahturin}). 

Let 
\bes
M_n=T_1\oplus T_2\oplus \ldots \oplus T_k\oplus\ldots 
\ees
be the standard degree grading of $M_n$, i.e., 
$T_k$ is the linear span of (non-associative) monomials in $x_1,x_2,\ldots,x_n$ of degree $k$. If 
\bes
f=f_s+f_{s+1}+\ldots+f_t \in M_n, f_j\in T_j, f_s\neq 0, f_t\neq 0, 
\ees
then $\deg(f)=t$ is the {\em degree} and $\mathrm{ldeg}(f)=s$ is the {\em lower degree} of $f$ in the variables $x_1,x_2,\ldots,x_n$.  Set also $\deg(0)=-\infty$ and $\mathrm{ldeg}(0)=\infty$. 

Let 
\bes
\phi=(x_1+f_1,x_2+f_2,\ldots,x_n+f_n), 
\ees
where $f_i\in M_n'$ for all $1\leq i\leq n$. Then we define $\deg(\phi)$ as the maximum of all $\deg(f_i), 1\leq i\leq n$, and $\mathrm{ldeg}(\phi)$ as the minimum of all $\mathrm{ldeg}(f_i), 1\leq i\leq n$. If $\deg(\phi)=\mathrm{ldeg}(\phi)=k\geq 2$ then we say that $\phi$ is {\em homogeneous of degree $k$}. 
In particular, if $k=2$ then $\phi$ is called {\em quadratic} and if $k=3$ then it is called {\em cubic}. 

Let $\mathrm{Aut}(M_n)$ be the group of all automorphisms of $M_n$ and let $\mathrm{TAut}(M_n)$ be the group of all tame automorphisms of $M_n$. Let also $\mathrm{ATAut}(M_n)$ be the group of all almost tame automorphisms of $M_n$. 
We have 
\bes
\mathrm{TAut}(M_n)\subseteq \mathrm{ATAut}(M_n)\subseteq \mathrm{Aut}(M_n). 
\ees

As usual, the group of all linear automorphisms of $M_n$ will be identified with the group of all invertible matrices $GL_n(K)$. Recall that the symmetric group $S_n$ is embedded into $GL_n(K)$ and any permutation $\pi\in S_n$ represents the automorphism $\pi=(x_{\pi(1)},\ldots,x_{\pi(n)})$. For example, the transposition $(23)$ represents the automorphism $(23)=(x_1,x_3,x_2,x_4,\ldots,x_n)$.

Obviously, the group of almost tame automorphisms $\mathrm{ATAut}(M_n)$ is generated by all linear automorphisms $GL_n(K)$ and all Chein automorphisms of the form  
\bee\label{f5}
(x_1+f,x_2,\ldots,x_n), \ f\in M_n'. 
\eee
Chein automorphisms of this form are described in \cite{AutoI}. 

The group of tame automorphisms $\mathrm{TAut}(M_n)$ is generated by all linear automorphisms $GL_n(K)$ and all automorphisms of the form (\ref{f5}) where $f\in M_n'$ does not depend on $x_1$.

The group of automorphisms $\mathrm{IAut}(M_n)$ of $M_n$ that induces the identity automorphism on $M_n/M_n'$ plays an important role in our study. Obviously, $\mathrm{IAut}(M_n)$ is the set of all automorphisms of $M_n$ of lower degree $\geq 2$. 
The elements of $\mathrm{IAut}(M_n)$ are called {\em $\mathrm{IA}$-automorphisms}.

Here we describe a set of generators of the group of tame $\mathrm{IA}$-automorphisms and almost tame $\mathrm{IA}$-automorphisms. 

\begin{lemma}\label{l1}
$(a)$ The group of tame $\mathrm{IA}$-automorphisms $\mathrm{TAut}(M_n)\cap \mathrm{IAut}(M_n)$ is generated by all linear conjugates of elementary automorphisms of the form (\ref{f5}) where $f\in M_n'$ does not depend on $x_1$;

$(b)$ The group almost tame $\mathrm{IA}$-automorphisms  $\mathrm{ATAut}(M_n)\cap \mathrm{IAut}(M_n)$ is generated by all linear conjugates of Chein automorphisms of the form (\ref{f5}). 
\end{lemma}
\Proof Every tame automorphism $\phi$ can be represented as
\bes
\phi=\a_1\phi_1\a_2\phi_2\ldots \phi_k\a_{k+1}, 
\ees
where $\phi_i$ is an elementary automorphism of the form (\ref{f5}) and $\a_i$ is linear for all $i$. Since $\phi\in \mathrm{IAut}(M_n)$  and $\phi_i\in \mathrm{IAut}(M_n)$ for all $i$ it follows that  
\bes
\a_1\a_2\ldots a_{k+1}=\mathrm{id}. 
\ees
 Then $\a_{k+1}=(\a_1\ldots \a_k)^{-1}$ and, consequently, 
\bes
\phi=(\a_1\phi_1\a_1^{-1})(\a_1\a_2\phi_2(\a_1\a_2)^{-1})\ldots (\a_1\ldots \a_k\phi_k(\a_1\ldots \a_k)^{-1})
\ees

The statement $(b)$ of the lemma can be proved similarly. $\Box$

 Let $Y_n$ be the abelian Lie algebra with a linear basis $y_1,y_2,\ldots,y_n$. The universal enveloping algebra $U(Y_n)$ of $Y_n$ is the polynomial algebra $U=U(Y_n)=K[y_1,\ldots,y_n]$ in the variables $y_1,y_2,\ldots,y_n$. Let 
\bes
U=U_0\oplus U_1\oplus \ldots \oplus U_k\oplus \ldots
\ees
be the degree grading of $U$, i.e., $U_k$ is the linear span of all monomials  of degree $k$ in $y_1,y_2,\ldots,y_n$.
If 
\bes
a=a_s+a_{s+1}+\ldots+a_t\in U, a_j\in U_j, a_s\neq 0, a_t\neq 0, 
\ees
then $\deg(a)=t$ is the {\em degree} and $\mathrm{ldeg}(a)=s$ is the {\em lower degree} of $a$ in the variables $y_1,y_2,\ldots,y_n$.
Set $\deg(0)=-\infty$ and $\mathrm{ldeg}(0)=\infty$ as usual.

Let $T_n$ be a free right $U$-module with basis $t_1,\ldots,t_n$. Turn the direct sum 
\bes
M=Y_n\oplus T_n 
\ees
  into a Lie algebra by 
\bes
[a+t,b+s]= tb-sa, 
\ees
where $a,b\in Y_n$ and $t, s\in T_n$. Then the subalgebra of $M$ generated by 
\bes
x_i=y_i+t_i, \ \ 1\leq i\leq n, 
\ees
is the free metabelian Lie algebra $M_n$ with the free generators $x_1,\ldots,x_n$. In the case of groups the corresponding embedding is called the Magnus embedding \cite{KSh,Magnus} and in the case of Lie algebras it follows from Shmel'kin's Theorem on wreath products \cite{Shmelkin73}. 

Notice that $M_n'\subseteq T_n$ and for any $m\in M_n'$ and $y_i$ we get 
\bes
m y_i=[m,x_i]. 
\ees

Every element $m\in M_n'$ can be written in the form 
\bes
m=\sum_{1\leq i<j\leq n} [x_i,x_j]a_{ij}, \quad a_{ij}\in U. 
\ees

Let $f\in M_n$ be an arbitrary element $M_n$. Then $f$ is uniquely represented as 
\bes
f=y+t, \quad y\in Y_n, \ t\in T_n. 
\ees
Moreover, $t$ is uniquely represented as 
\bes
t=t_1d_1+\ldots+t_nd_n, \ \ \ \ d_1,\ldots,d_n\in K[y_1,\ldots,y_n]=U. 
\ees
Set  $d_i=\frac{\partial f}{\partial x_i}$ for all $1\leq i\leq n$ and 
\bes
\partial(f)=(d_1,\ldots,d_n)^t=\big(\frac{\partial f}{\partial x_1},\ldots,\frac{\partial f}{\partial x_n}\big)^t, 
\ees
where $t$ means the transpose. We call $\frac{\partial f}{\partial x_1},\ldots,\frac{\partial f}{\partial x_n}$ the {\em Fox derivatives of $f$}  (see \cite{Shpilrain,U93}). 

Let $e_1,e_2,\ldots,e_n$ be the columns of the identity matrix $I=I_n$ of order $n$. 

The statement of the next lemma is well known (see \cite[Lemma 5]{U93},\cite{95SMJ}). 
\begin{lemma}\label{l2} \cite[Lemma 1]{AutoI} Let $a=(a_1,a_2,\ldots,a_n)^t\in (U^n)^t$ be an arbitrary column of dimension $n$ over $U$. 
Then the following conditions are equivalent. 

$(1)$ $a=\d(f)$ for some $f\in M_n'$; 

$(2)$  $Ya=y_1a_1+\ldots +y_na_n=0$ where $Y=(y_1,\ldots,y_n)$; 

$(3)$ $a=\sum_{i<j} (e_iy_j-e_jy_i)a_{ij}, \ a_{ij}\in U$. 
\end{lemma}

Every endomorphism $\phi$ of $M_n$ induces an endomorphism of $Y_n=M_n/M_n'$, which can be uniquely extended to an endomorphism of $U=U(Y_n)$. This endomorphism will be denoted again by $\phi$ for simplicity of notation. If $f\in U$, then the image $\phi(f)$ will be denoted by $f^{\phi}$. 

We define the Jacobian matrix of $\phi=(f_1,\ldots,f_n)$  by 
\bes
J(\phi)=\left[\frac{\partial f_j}{\partial x_i}\right]_{1\leq i,j\leq n}=\left[\partial(f_1) \  \ldots \ \partial(f_n)\right].
\ees
If $\phi=(f_1,\ldots,f_n)$ and $\psi=(g_1,\ldots,g_n)$ then 
\bes
\phi\circ\psi(x_i)=g_i(f_1,\ldots,f_n)
\ees
for all $1\leq i\leq n$. The corresponding Chain Rule gives \cite{Shpilrain,U93,95SMJ} that 
\bee\label{f6}
J(\phi\circ\psi)=J(\phi) J(\psi)^{\phi}. 
\eee

The statement of the next lemma is also well known (see \cite{Shpilrain},\cite[Theorem 2]{U93}). 
\begin{lemma}\label{l3} \cite[Lemma 2]{AutoI}
An endomorphism $\phi$ of $M_n$ is an automorphism if and only if the Jacobian matrix $J(\phi)$ is invertible.  
\end{lemma}

Here is a description of Chein automorphisms. 

\begin{lemma}\label{l4} \cite[Lemma 3]{AutoI} An endomorphism of $M_n$ of the form (\ref{f5}) 
 is an automorphism if and only if $f$ belongs to the ideal generated by all $[x_s,x_t]$, where $s,t\neq 1$.   
\end{lemma}

Consider the set of Chein automorphisms of the form 
\bee\label{f7}
(x_1+[x_2,x_3]a, x_2,\ldots,x_n), \ a\in U. 
\eee
These automorphisms are elementary if $a$ does not depend on $y_1$. 

In this paper we denote by $\mathrm{Mat}_n(U)$ the set of all square matrices of order $n$ over $U$ (since the more standard notation $M_n(U)$ can be confused with the notation of the free metabelian Lie algebra $M_n$).
As usual, denote by $GL_n(U)$ and $SL_n(U)$ the general linear group and the special linear group of matrices of order $n$ over $U$, respectively. 

We mentioned above that the group $GL_n(K)$ will be identified with the group of all linear automorphisms of $M_n$. We meant that each linear automorphism $\a$ will be identified with its Jacobian matrix $J(\a)$. This means that $\a$ is a matrix and instead of $J(\a)$ in our calculations we will write $\a$.

\begin{lemma}\label{l5} $(a)$ If $\phi, \psi\in  \mathrm{IAut}(M_n)$ then 
\bes
J(\phi\psi)=J(\phi) J(\psi); 
\ees

$(b)$  If $\phi\in  \mathrm{IAut}(M_n)$ and $\a\in GL_n(K)$ then 
\bes
J(\a \phi \a^{-1})=\a J(\phi)^{\a} \a^{-1}; 
\ees

$(c)$ 
 The map 
\bes
J : \mathrm{IAut}(M_n) \to SL_n(U)
\ees
is an injective homomorphism of groups and $G=J(\mathrm{IAut}(M_n))$ is the set of all matrics from  $SL_n(U)$ of the form 
\bes
I+A, 
\ees
where $A\in \mathrm{Mat}_n(U)$ is a matrix such that $YA=0$. 
\end{lemma}
\Proof Since $\phi$ induces an identity automorphism on $U$, $(a)$ is a direct consequence of (\ref{f6}). Using (\ref{f6}) we also get
\begin{align*}
J(\a \phi \a^{-1})=J(\a (\phi \a^{-1}))=J(\a)J(\phi \a^{-1})^{\a}=J(\a)(J(\phi) J(\a^{-1})^{\phi})^{\a}
=J(\a)(J(\phi) J(\a^{-1}))^{\a}\\
=J(\a)J(\phi)^{\a} J(\a^{-1})^{\a}=J(\a)J(\phi)^{\a} J(\a^{-1})=J(\a)J(\phi)^{\a} J(\a)^{-1}=\a J(\phi)^{\a} \a^{-1}. 
\end{align*}
In this calculation we again used that $\phi$ induces an identity automorphism on $U$, $J(\a^{-1})\in GL_n(K)$ is a constant matrix, and that $J(\a^{ -1})=J(\a)^{-1}$ by (\ref{f6}).

 If $\phi\in \mathrm{IAut}(M_n)$, then 
\bes
\phi=(x_1+f_1,x_2+f_2,\ldots,x_n+f_n), 
\ees
 where $f_i\in M_n'$ 
 for all $i$, and 
\bes
J(\phi)=I+A, \ A=\left[\partial(f_1) \  \ldots \ \partial(f_n)\right].
\ees
We get $YA=0$ by Lemma \ref{l2}. If $J(\phi)=I$ then $\d(f_i)=0$ for all $i$. Notice that for any $f\in M_n'$ we get $f=(t_1,\ldots,t_n)\d(f)$ by the definition of $\d(f)$. Consequently, $f_i=(t_1,\ldots,t_n)\d(f_i)=0$. 
 This means $\phi=\mathrm{id}$, i.e., $J$ is injective. 

Conversely, let $A\in \mathrm{Mat}_n(U)$ such that $YA=0$. By Lemma \ref{l2} there exist $g_1,g_2,\ldots,g_n\in M_n'$ such that $A=\left[\partial(g_1) \  \ldots \ \partial(g_n)\right]$. 
Set $\psi=(x_1+g_1,x_2+g_2,\ldots,x_n+g_n)$. Then $J(\psi)=I+A$. If $I+A\in SL_n(U)$ then $\psi$ is an automorphism by Lemma \ref{l3} and $\psi\in \mathrm{IAut}(M_n)$. $\Box$

\section{Matrix notations}

\hspace*{\parindent} 

Using Lemma \ref{l5}, we will represent automorphisms from $\mathrm{IAut}(M_n)$ by their Jacobian matrices and make calculations in the language of Jacobian matrices. The results of \cite{AutoI}, originally formulated in the language of automorphisms, will be reformulated in the language of matrices, while the results of the present paper, originally formulated in the language of matrices, will also be reformulated in the language of automorphisms.

Thus $G$ represents the set of all Jacobian matrices of automorphisms from $\mathrm{IAut}(M_n)$ as in Lemma \ref{l5}. 
Set $T=J( \mathrm{TAut}(M_n)\cap \mathrm{IAut}(M_n))$ and $S=J(\mathrm{ATAut}(M_n)\cap \mathrm{IAut}(M_n))$. 
We have 
\bes
T\subseteq S \subseteq G. 
\ees 
In order to prove that $ \mathrm{ATAut}(M_n)=\mathrm{Aut}(M_n)$ it is sufficient to prove that $ \mathrm{ATAut}(M_n)\cap \mathrm{IAut}(M_n)=\mathrm{IAut}(M_n)$ since both groups include all linear automorphisms. Thus, it is sufficient to prove that $S=G$.

Let $E_{ij}$ denote the square matrix with $1$ in the $(i,j)$-entry and $0$ elsewhere. The order of $E_{ij}$ is usually clear from the context. Recall that 
\bes
E_{ij}E_{kl}=\delta_{jk}E_{il}, 
\ees
where $\delta_{jk}$ is the Kronecker delta. 
A matrix of the form 
\bes
I+aE_{ij}, \ i\neq j, \ a\in U,
\ees
is called {\em elementary}. 

If $i<j$, then for any $k$ denote by $F_{ijk}$ the $n\times n$ matrix with $y_j$ in position $(i,k)$, $-y_i$ in position $(j,k)$, and zeros in other places, that is,
\bes
F_{ijk}=y_jE_{ik}-y_iE_{jk}. 
\ees
 In other words, the $k$th column of $F_{ijk}$ is equal to $\d([x_i,x_j])=e_iy_j-e_jy_i$, and the rest are zeros. It is easy to check that 
\bes
F_{ijk}F_{stm}=\left\{\begin{array}{ccc}
y_tF_{ijm} & {\rm if} & k=s,\\
-y_sF_{ijm}      & {\rm if} & k=t,\\
0     & {\rm if} & k\neq s,t\\
\end{array}\right., 
\ees
\bes
E_{ij}F_{kst}=y_sE_{ij}E_{kt}-y_kE_{ij}E_{st}=
\left\{\begin{array}{ccc}
y_sE_{it} & {\rm if} & j=k,\\
-y_kE_{it}      & {\rm if} & j=s,\\
0     & {\rm if} & j\neq k,s\\
\end{array}\right., 
\ees
and 
\bes
F_{kst}E_{ij}=y_sE_{kt}E_{ij}-y_kE_{st}E_{ij}=
\left\{\begin{array}{ccc}
F_{ksj} & {\rm if} & t=i,\\
0     & {\rm if} & t\neq i\\
\end{array}\right..
\ees

We refer to these three formulas as the {\em multiplication tables} but often use them without reference. 
Properties of the product of matrices of types 
\bes
I+A+B=(I+A)(I+B),
\ees
which is true if $AB=0$, and 
\bes
I+(a+b)C=(I+aC)(I+bC), 
\ees
which is valid if $C^2=0$, to simplify terminology, will be called {\em additivity}.  
For example, if 
\bes 
f=\sum_{1<i<j}[x_i,x_j]a_{ij},
\ees
 then the Jacobian matrix of the endomorphism (\ref{f5}) is
\bes
I +\sum_{i<j}a_{ij} F_{ij1}, 
\ees
and, by additivity, 
\bes
I +\sum_{i<j}a_{ij} F_{ij1}=\prod_{i<j}(I+a_{ij} F_{ij1})
\ees
since 
\bes
F_{ij1}F_{st1}=0 
\ees
if $i,j,s,t\neq 1$. 
Moreover, if $1<i<j$, then 
\bes
I +(a+b)F_{ij1}=(I +aF_{ij1})(I +bF_{ij1}). 
\ees

Consider the set of matrices of the form 
\bee\label{f8}
I +aF_{231},  \ a\in U. 
\eee
These matrices are the Jacobian matrices of the Chein automorphisms defined in (\ref{f7}).

Here are two remarks:

$(a)$ If $A\in G$ and $\a\in GL_n(K)$ then $\a A^{\a}\a^{-1}\in G$ by Lemma \ref{l5}$(b)$. This matrix will be called the {\em linear conjugate of $A$ by $\a$}. In many places, conjugation by linear automorphisms helps to ``linearize'' automorphisms \cite{KBKA}.

$(b)$  If $J(\phi)=I+A\in S$ and $A$ is a homogeneous matrix, then, usually, it is easy to show that $I+\mu A\in S$ for any $\mu\in K$. This can be done by adding some coefficients in the proofs. For example, this can be done when conjugating matrices (\ref{f8}) using $(x_1,\mu x_2,x_3,\ldots,x_n)$. For this reason we ignore the difference between the coefficients of this type.

\begin{lemma}\label{l6} $(a)$ 
The group $T$ is generated by all linear conjugates of matrices of the form (\ref{f8}) where $a$ does not depend on $y_1$. 

$(b)$ The group $S$ is generated by all linear conjugates of matrices of the form (\ref{f8}). 
\end{lemma}
\Proof By Lemma \ref{l1}, every automorphism $\mathrm{TAut}(M_n)\cap \mathrm{IAut}(M_n)$ is a product of linear conjugates of automorphisms of the form (\ref{f5}) where $f$ does not depend on $x_1$. Observe that every such automorphism is itself a product of linear conjugates of automorphisms of the form (\ref{f7}) with $a\in K[y_2,\ldots,y_n]$. If $\t$ is an automorphism of the form (\ref{f7}) then $J(\t)$ is a matrix of the form (\ref{f8}).  If $\a\in GL_n(K)$ then 
\bes
J(\a\t\a^{-1})=\a J(\t)^{\a}\a^{-1}
\ees
by Lemma \ref{l5}, that is, $J(\a\t\a^{-1})$ is the linear conjugate of $J(\t)$ by $\a$. If $\phi,\psi\in \mathrm{IAut}(M_n)$ then $J(\phi\psi)=J(\phi)J(\psi)$ 
 by the same lemma. This proves $(a)$. 

The statement $(b)$ can be checked similarly. $\Box$

Let's turn some results of \cite{AutoI} into matrix form. 

\begin{proposition}\label{p1} Let $n\geq 4$ and $K$ is a field of characteristic $\neq 3$ or $n\geq 5$. 
	
	$(1)$ Every matrix of the form (\ref{f8}) belongs to $T$ if $\mathrm{ldeg}(a)\geq 2$; 
	
	$(2)$ Every matrix of the form 
	\bes
	I+ay_1F_{121}+ay_2F_{122}
	\ees
	belongs to $T$ if $\mathrm{ldeg}(a)\geq 2$.
\end{proposition}
\Proof It is clear that the Jacobian matrix of a Chein automorphism of the form (\ref{f7}) has the form (\ref{f8}). By Theorem 1 of \cite{AutoI}, every automorphism of the form (\ref{f7}) is tame whenever $\mathrm{ldeg}(a)\geq 2$. Therefore, the statement $(1)$ of the lemma follows.

By \cite[Lemma 4]{AutoI}, all automorphisms of the form 
\bes
(x_1+[x_1,x_2] y_1a, x_2+[x_1,x_2] y_2a,x_3,\ldots,x_n)
\ees
are tame if $\mathrm{ldeg}(a)\geq 2$. Note that statement $(2)$ of the lemma is the matrix formulation of this statement. $\Box$

\begin{corollary}\label{c1} \cite[Corollary 1]{AutoI} Let $n\geq 4$ and $K$ is a field of characteristic $\neq 3$ or $n\geq 5$.
	The group of all almost tame automorphisms $\mathrm{ATAut}(M_n)$ is generated, modulo all linear automorphisms, by the quadratic automorphism (\ref{f1}) and the  cubic automorphism (\ref{f2}).
\end{corollary}

\begin{proposition}\label{p2} Every matrix of the form 
	\bes
	I+ay_1F_{121}+ay_2F_{122}
	\ees
	with $\mathrm{ldeg}(a)\geq 1$ belongs to $S$.  
\end{proposition}
\Proof By \cite[Lemma 5]{AutoI}, all automorphisms of the form 
\bes
(x_1+[x_1,x_2] y_1a, x_2+[x_1,x_2] y_2a,x_3,\ldots,x_n)
\ees
are almost tame if $\mathrm{ldeg}(a)\geq 1$. The matrix formulation of this statement is precisely the statement of the proposition. $\Box$

Observe that Proposition \ref{p2} does not include only the matrix 
\bee\label{f9}
I +y_1F_{121}+y_2F_{122}
\eee
of the same type, which is the Jacobian matrix of (\ref{f3}). We consider this matrix separately in Section 6.

\section{Special matrices}

\hspace*{\parindent}

We use matrices instead of automorphisms in many places because it makes the calculations easier to understand. To some extent this is due to the application of Suslin's theorem \cite{Suslin}. 

The main aim of this section is to prove that every matrix of the form
\bee\label{f10}
I +\begin{bmatrix}
y_nC & \begin{matrix}
0 \\
. \\
.\\
.\\
0
\end{matrix} \\
\begin{matrix}
* & * & \ldots & *
\end{matrix} & 0
\end{bmatrix}\in G, 
\eee
where $C\in \mathrm{Mat}_{n-1}(U)$, 
 belongs to $S$.   We will provide several lemmas to the proof of this result. A matrix of the form (\ref{f10}) will be called {\em special}.

If $J$ is an ideal of $U$ then the natural map $U\to U/J$ induces the natural map $SL_k(U)\to SL_k(U/J)$. Denote by  $SL_k(U,J)$ the kernel of $SL_k(U)\to SL_k(U/J)$.

\begin{lemma}\label{l7}
Let $B$ be any matrix from $SL_{n-1}(U,y_nU)$. Then $B=I_{n-1}+y_nC$ for some matrix $C\in \mathrm{Mat}_{n-1}(U)$ and there exists a unique special matrix $Z$ such that 
\bes
Z=I+\left[\begin{array}{cc}
 y_nC  & 0  \\
 *  &  0  \\
 \end{array}\right].
\ees
\end{lemma}
\Proof Obviously, $B=I_{n-1}+y_nC$ for some matrix $C\in\mathrm{Mat}_{n-1}(U)$. If $C=[c_{ij}]$ then 
\bes
Z=I +\begin{bmatrix}
y_nC & \begin{matrix}
0 \\
. \\
.\\
.\\
0
\end{matrix} \\
\begin{matrix}
-\sum_{i=1}^{n-1} y_ic_{i1} & -\sum_{i=1}^{n-1} y_ic_{i2} & \ldots & -\sum_{i=1}^{n-1} y_ic_{i (n-1)}
\end{matrix} & 0
\end{bmatrix}
\ees
belongs to $G$ by Lemma \ref{l2}. Lemma \ref{l2} also implies that the last row of $Z$ is uniquely determined by $C$. This implies the uniqueness of $Z$. $\Box$

 If $B\in SL_n(U)$ is a matrix of the form $B=
\begin{bmatrix}
C & 0 \\
0 & I_{n-k}
\end{bmatrix} $, where $C\in SL_k(U)$ for some $1\leq k\leq n$, then we say that $B$ is {\em $k$-reducible} and 
$B_k=C$ will be called its {\em $k$-reduction}. Calculations with $k$-reducible matrices from $G$ can be done by using their $k$-reductions. 

Set $Z_{n-1}=B$ for matrices $Z,B$ satisfying the conditions of Lemma \ref{l7}. We use this notation for special matrices despite the fact that $Z$ is not $n-1$-reducible in general. But Lemma \ref{l7} allows us to use the matrix $B$ instead of $Z$ in some calculations. Moreover, if $B$ is $k$-reducible for some $1\leq k\leq n-1$, then $B_k$ will also be called the {\em $k$-reduction of the special matrix $Z$} and the special matrix $Z$ will be called {\em special $k$-reducible}.

 Observe that if $A$ and $B$ are two special $k$-reducible matrices, then $AB$ is a special $k$-reducible and $(AB)_k=A_kB_k$ by Lemma \ref{l5}$(a)$ and Lemma \ref{l7}.

\begin{proposition}\label{p3} Every matrix of the form   
\bes
A(h,f,g)=I +hfgF_{1n1}+hg^2F_{2n1}-hf^2F_{1n2}-hfgF_{2n2}
\ees
belongs to $S$. 
\end{proposition}
\Proof Theorem 3 in \cite{AutoI} says that every automorphism of the form 
\bes
(x_1+[x_1,x_n]hfg+[x_2,x_n]hg^2, x_2-[x_1,x_n]hf^2-[x_2,x_n]hfg,x_3,\ldots,x_n)
\ees
is almost tame. The matrix formulation of this result is precisely the statement of the proposition. 
 $\Box$

Let us recall some results of Cohn \cite{Cohn66} and Suslin \cite{Suslin} on generators of special linear groups over polynomial algebras. 

In 1966 P. Cohn \cite{Cohn66} proved that the matrix  
\bes
\begin{bmatrix}
    1+x_1x_2 & x_2^2\\
-x_1^2 & 1-x_1x_2\\
\end{bmatrix}
\ees
is not a product of elementary matrices over the polynomial algebra $K[x_1,x_2]$. 

Denote by $E_k(U)$ the subgroup of $SL_k(U)$ generated by all elementary matrices, that is, by all matrices of the form 
\bes
I_k+aE_{ij}, \ a\in U. 
\ees
 The well-known Suslin's Theorem \cite{Suslin} states that if $k\geq 3$ then $SL_k(U)=E_k(U)$, that is, every matrix from  $SL_k(U)$ is a product of elementary matrices. If $J$ is an ideal of $U$ then denote by $E_k(U,J)$ the normal closure in $E_k(U)$ of the set of all elementary matrices of the form $I_k+aE_{ij}$, where $i\neq j$ and $a\in J$. 

The next corollary of Suslin's results is formulated in \cite{BM85}. 
\begin{proposition}\label{p4}
Let  $k\geq 3$ and $J=y_nU$. Then 

(1) $SL_k(U,J)=E_k(U,J)$;

(2) $E_k(U,J)$  is generated by the set of all matrices of the form 
\bes
I_k+h(fe_i+ge_j)(ge_i-fe_j)^t=I_k+hfg(E_{ii}-E_{jj})+hg^2E_{ji}-hf^2E_{ij}
\ees
where $i<j$, $h\in J$, and $f,g\in U$.
\end{proposition}
\Proof We repeat a short proof given in \cite{BM85}. By \cite{Suslin}, every matrix $\Phi\in SL_k(U,J)$ is a product of elementary matrices $\Phi=E_1\ldots E_s$. Any $a\in U=K[y_1,\ldots,y_n]$ can be written in the form 
\bes
a=y_nb+c, \ \ b\in U, c\in K[y_1,\ldots,y_{n-1}], 
\ees
and, consequently, 
\bes
I_k+aE_{ij}=(I_k+y_nbE_{ij})(I_k+cE_{ij}). 
\ees
Thus, $E_i=H_iG_i$ where $H_i\in E_k(U,J)$ and $G_i\in E_k(K[y_1,\ldots,y_{n-1}])$. Then
\bes
\Phi=H_1(G_1H_2G_1^{-1})(G_1G_2H_3G_2^{-1}G_1^{-1})\ldots(G_1\ldots G_{s-1}H_sG_{s-1}^{-1}\ldots G_1^{-1})(G_1\ldots G_s). 
\ees
Substituting $y_n=0$ we get $G_1\ldots G_s=1$. Since 
\bes
G_1\ldots G_{i-1}H_iG_{i-1}^{-1}\ldots G_1^{-1}\in E_k(U,J)
\ees
it follows that $\Phi\in E_k(U,J)$.

By \cite[Corollary 1.4]{Suslin}, the group $E_k(U,J)$ is generated by the set of all matrices of the form 
\bes
I_k+hu(u_je_i-u_ie_j)^t
\ees
where $h\in J$, $u=(u_1,\ldots,u_k)$ is a unimodular column over $U$, and $1\leq i<j\leq k$. Recall that a column $u$ is called {\em unimodular} if the ideal generated by its components contains the identity element. We have 
\begin{align*}
I_k+hu(u_je_i-u_ie_j)^t=(I_k+h(u_ie_i+u_je_j)(u_je_i-u_ie_j)^t)\\
\times \prod_{l\neq i,j} (I_k+hu_le_l(u_je_i)^t)
\prod_{l\neq i,j} (I_k-hu_le_l(u_ie_j)^t).
\end{align*}
This proves the statement $(2)$ of the proposition. 
$\Box$

Observe that if $i=1,j=2$, then the matrix from Proposition \ref{p4}$(2)$ is $2$-reducible and its $2$-reduction is 
\bes
I_2+h\left[\begin{array}{cc}
 fg  & -f^2 \\
 g^2 &  -fg  \\
 \end{array}\right], 
\ees
which is a generalization of Cohn's matrix mentioned above.

\begin{proposition}\label{p5} 
Every special matrix belongs to $S$. 
\end{proposition}
\Proof Since $n\geq 4$ it follows by Proposition \ref{p4} that $SL_{n-1}(U,y_nU)$ is generated by the set of all matrices
\bes
B=I_{n-1}+y_nh(fe_i+ge_j)(ge_i-fe_j)^t, \ 1\leq i<j\leq n-1. 
\ees
It is sufficient to prove that for any such matrix $B$ there exists a special matrix $A$ from $S$ such that its special $(n-1)$-reduction is $B$.  
Without loss of generality we may assume that $i=1$ and $j=2$, that is, 
\bes
B=I_{n-1}+y_nh(fe_1+ge_2)(ge_1-fe_2)^t. 
\ees
Then for 
\bes
A=I +hfgF_{1n1}+hg^2F_{2n1}-hf^2F_{1n2}-hfgF_{2n2} 
\ees
we get $A_{n-1}=B$ and $A$ belongs to $S$ by Proposition \ref{p3}. $\Box$

The proposition can be reformulated in the language of automorphisms as follows.

\begin{corollary}\label{c2}
Every automorphism of the form 
\bes
(x_1+f_1,x_2+f_2,\ldots,x_{n-1}+f_{n-1}, x_n), 
\ees
where $f_1,f_2,\ldots,f_{n-1}$ belong to the ideal generated by $[x_1,x_n],\ldots,[x_{n-1},x_n]$, is almost tame.  
\end{corollary}

\section{Inductive lemma}

\hspace*{\parindent} 

Our aim is to reduce the number of variables in a matrix $A\in G$ by using matrices from $S$. More exactly, for any matrix $A\in G$ we try to find a matrix $B\in S$ such that $B^{-1}A$ depends on a less number of variables than $A$. 
The following lemma reduces the number of variables in a particular case. The proof follows that of the corresponding statement for groups \cite[pp.~99--100]{BM85}. Unfortunately, we cannot establish it in its full generality at this stage, since exponential automorphisms of degree 3 are not yet available. We therefore begin by proving a slightly restricted version. In the next section, we show that exponential automorphisms of degree 3 are almost tame, which enables us to prove the lemma in full generality. Fortunately, the proof presented below applies without modification to the general case. 

\begin{lemma}\label{l8}
Let $B=I+[a_{ij}]\in SL_n(U,\sum_{p=s}^ny_pU)\cap G$ be such that   
\bes
a_{in}\in
\sum_{p=s}^{n-1}y_pU+y_n^2U
\ees
 and  
 \bes
 \mathrm{ldeg}(\overline{a}_{in})\geq 3
 \ees
  for all $1\leq i\leq n-1$, where $\overline{a}_{in}$ is obtained from $a_{in}$ by setting $y_s=y_{s+1}=\ldots =y_{n-1}=0$. Then there exists $C\in S$ such that $C^{-1}B\in SL_n(U,\sum_{p=s}^{n-1}y_pU)$. 
\end{lemma}
\Proof For any $a\in U$ and any matrix $R$ over $U$ let $\overline{a}$ and $\overline{R}$ be the result of setting $y_s=y_{s+1}=\ldots =y_{n-1}=0$ in $a$ and $R$ respectively. 

Set $A=[a_{ij}]$. Recall that $(y_1,\ldots,y_n)[a_{ij}]=0$ by Lemma \ref{l2}.

Let 
\bes
M=
\begin{bmatrix}
I_{n-1} & \begin{matrix}
0 \\
. \\
.\\
.\\
0
\end{matrix} \\
\begin{matrix}
y_1 & y_2 & \ldots & y_{n-1}
\end{matrix} & y_n 
\end{bmatrix}. 
\ees
Then
\bes
\overline{M}=
\begin{bmatrix}
I_{n-1} & \begin{matrix}
0 \\
. \\
.\\
.\\
0
\end{matrix} \\
\begin{matrix}
y_1 & \ldots & y_{s-1} & 0 & \ldots & 0
\end{matrix} & y_n 
\end{bmatrix}. 
\ees
Observe that 
\bes
M^{-1}=
\begin{bmatrix}
	I_{n-1} & \begin{matrix}
		0 \\
		. \\
		.\\
		.\\
		0
	\end{matrix} \\
	\begin{matrix}
		-y_1y_n^{-1} & -y_2y_n^{-1} & \ldots & -y_{n-1}y_n^{-1}
	\end{matrix} & y_n^{-1} 
\end{bmatrix}. 
\ees

 We have 
\begin{align*}
M B M^{-1}=I_n+ \begin{bmatrix}
[a_{ij}-y_n^{-1}y_ja_{in}] & \begin{matrix}
y_n^{-1}a_{1n} \\
. \\
.\\
.\\
y_n^{-1}a_{(n-1) n}
\end{matrix} \\
\begin{matrix}
0 & \ldots & 0
\end{matrix} & 0
\end{bmatrix}\\
=\begin{bmatrix}
I_{n-1} & \begin{matrix}
y_n^{-1}a_{1n} \\
. \\
.\\
.\\
y_n^{-1}a_{(n-1) n}
\end{matrix} \\
\begin{matrix}
0 & \ldots & 0
\end{matrix} & 1
\end{bmatrix}
\left(I_n+ \begin{bmatrix}
[a_{ij}-y_n^{-1}y_ja_{in}] & \begin{matrix}
0 \\
. \\
.\\
.\\
0
\end{matrix} \\
\begin{matrix}
0 & \ldots & 0
\end{matrix} & 0
\end{bmatrix}\right)
\end{align*}
and
\bes
\overline{M} \ \overline{B} \ \overline{M}^{-1}
=\begin{bmatrix}
I_{n-1} & \begin{matrix}
y_n^{-1}\overline{a}_{1n} \\
. \\
.\\
.\\
y_n^{-1}\overline{a}_{(n-1) n}
\end{matrix} \\
\begin{matrix}
0 & \ldots & 0
\end{matrix} & 1
\end{bmatrix}
\left(I_n+ \begin{bmatrix}
[\overline{a}_{ij}-y_n^{-1}\overline{y_j} \, \overline{a}_{in}] & \begin{matrix}
0 \\
. \\
.\\
.\\
0
\end{matrix} \\
\begin{matrix}
0 & \ldots & 0
\end{matrix} & 0
\end{bmatrix}\right). 
\ees
Then
\bes
\overline{M} \ \overline{B} \ \overline{M}^{-1}
=\begin{bmatrix}
I_{n-1} & \begin{matrix}
y_n\overline{d}_{1n} \\
. \\
.\\
.\\
y_n\overline{d}_{(n-1) n}
\end{matrix} \\
\begin{matrix}
0 & \ldots & 0
\end{matrix} & 1
\end{bmatrix}
\left(I_n+ \begin{bmatrix}
[y_n\overline{d}_{ij}] & \begin{matrix}
0 \\
. \\
.\\
.\\
0
\end{matrix} \\
\begin{matrix}
0 & \ldots & 0
\end{matrix} & 0
\end{bmatrix}\right),
\ees
where $\overline{d_{ij}}$ are defined by 
\bes
\overline{a_{in}}=y_n^2\overline{d_{in}}, \quad  1\leq i\leq n-1, 
\ees
and 
\bes
\overline{a_{ij}}-y_n^{-1}\overline{y_j}\overline{a_{in}}=y_n\overline{d_{ij}}, \quad 1\leq i,j\leq n-1
\ees
Such elements $\overline{d}_{ij}\in K[y_1,\ldots,y_{s-1},y_n]$ exist because, by the assumption of the lemma, 
\bes
\overline{a}_{in}\in y_n^2K[y_1,\ldots,y_{s-1},y_n] \quad 1\leq i\leq n-1, 
\ees 
and 
\bes
\overline{a}_{ij}\in y_nK[y_1,\ldots,y_{s-1},y_n] \quad 1\leq i,j\leq n-1. 
\ees 

Moreover,  $\mathrm{ldeg}(\overline{d}_{kn})\geq 1$ for all $1\leq k\leq n-1$ since $\mathrm{ldeg}(\overline{a}_{kn})\geq 3$. For any $1\leq k\leq n-1$, let 
\begin{align*}
\delta_k=I_n+\begin{bmatrix}
0 & 0 & \ldots & 0 \\
. & . & \ldots & . \\
. & . & \ldots & . \\
. & . & \ldots & . \\
0 & 0 & \ldots & 0 \\
y_1y_n\overline{d}_{kn} & y_2y_n\overline{d}_{kn} & \ldots & y_n^2\overline{d}_{kn} \\
0 & 0 & \ldots & 0 \\
. & . & \ldots & . \\
. & . & \ldots & . \\
. & . & \ldots & . \\
0 & 0 & \ldots & 0 \\
-y_1y_k\overline{d}_{kn} & -y_2y_k\overline{d}_{kn} & \ldots & -y_ny_k\overline{d}_{kn}
\end{bmatrix}\\
=\prod_{l\neq k,n} (I_n+y_l\overline{d}_{kn} F_{knl})\cdot (I_n+y_k\overline{d}_{kn}F_{knk}+y_n\overline{d}_{kn}F_{knn}). 
\end{align*}
Note that all factors of the product $\prod_{l\neq k,n} (I_n+y_l\overline{d}_{kn} F_{knl})$ belong to $S$ and 
\bes
I_n+y_k\overline{d}_{kn}F_{knk}+y_n\overline{d}_{kn}F_{knn}\in S
\ees
by Proposition \ref{p2} since $\mathrm{ldeg}(\overline{d}_{kn})\geq 1$ (this is the only place where the restriction $\mathrm{ldeg}(\overline{a}_{kn})\geq 3$ was used). Consequently, $\delta_k\in S$. We have 
\begin{align*}
	M\delta_k M^{-1}=I_n+\begin{bmatrix}
		0 & 0 & \ldots & 0 \\
		. & . & \ldots & . \\
		. & . & \ldots & . \\
		. & . & \ldots & . \\
		0 & 0 & \ldots & 0 \\
		0 & 0 & \ldots & y_n\overline{d}_{kn} \\
		0 & 0 & \ldots & 0 \\
		. & . & \ldots & . \\
		. & . & \ldots & . \\
		. & . & \ldots & . \\
		0 & 0 & \ldots & 0 \\
		0 & 0 & \ldots & 0 \\
	\end{bmatrix}.
\end{align*}

Therefore  
\bes
\overline{M} \overline{\delta}_1\ldots \overline{\delta}_{n-1} \overline{M}^{-1}=
\begin{bmatrix}
I_{n-1} & \begin{matrix}
y_n\overline{d}_{1n} \\
. \\
.\\
.\\
y_n\overline{d}_{(n-1) n}
\end{matrix} \\
\begin{matrix}
0 &  \ldots & 0
\end{matrix} & 1 
\end{bmatrix}
\ees
and 
\bes
\overline{M} (\overline{\delta}_1\ldots \overline{\delta}_{n-1})^{-1}\overline{B} \ \overline{M}^{-1}= I_n+
\begin{bmatrix}
[y_n\overline{d}_{ij}] & \begin{matrix}
0 \\
. \\
.\\
.
\end{matrix} \\
\begin{matrix}
0 &  \ldots 
\end{matrix} & 0 
\end{bmatrix}. 
\ees
Let 
\bes
\epsilon=I_n+\begin{bmatrix}
y_n\overline{d}_{11} & \ldots & y_n\overline{d}_{1 (n-1)} & 0 \\
. &  \ldots & . & . \\
. &  \ldots & . & . \\
. &  \ldots & . & . \\
y_n\overline{d}_{(n-1) 1} & \ldots & y_n\overline{d}_{(n-1) (n-1)} & 0 \\
-\sum_{j=1}^{n-1}y_j\overline{d}_{j 1} & \ldots & -\sum_{j=1}^{n-1}y_j\overline{d}_{j (n-1)} & 0
\end{bmatrix}. 
\ees
By Proposition \ref{p5}, $\epsilon\in S$, and 
\bes
\overline{M\epsilon M^{-1}}=\overline{M} \overline{\epsilon} \overline{M}^{-1}=I_n+\begin{bmatrix}
[y_n\overline{d}_{ij}] & \begin{matrix}
0 \\
. \\
.\\
.
\end{matrix} \\
\begin{matrix}
0 &  \ldots 
\end{matrix} & 0 
\end{bmatrix}. 
\ees
Set $C=\delta_1\ldots \delta_{n-1}\epsilon$. Then  
\bes
\overline{M} \ \overline{C}^{-1} \ \overline{B} \ \overline{M}^{-1}=\overline{M} \ \overline{\epsilon}^{-1} \ (\overline{\delta}_1\ldots \overline{\delta}_{n-1})^{-1}\overline{B} \ \overline{M}^{-1}=I_n. 
\ees
Consequently, $\overline{C}^{-1} \ \overline{B}=I_n$. This means that $C^{-1}B\in SL_n(U,\sum_{p=s}^{n-1}y_pU)$.   $\Box$

\section{Exceptional automorphism}

\hspace*{\parindent}

This section is devoted to proving that the matrix (\ref{f9}) belongs to $S$, or equivalently, that the automorphism (\ref{f3}) is almost tame. Originally, the proof was given in the language of matrices. In response to the referee's comments, we present the same calculations again in the language of automorphisms.  

In 1983, Anick \cite{Anick} proved that every automorphism (even every endomorphism with Jacobian determinant $1$) of the polynomial algebra $K[x_1,\ldots,x_n]$ over a field $K$ of characteristic zero can be approximated by tame automorphisms in the power series topology. This can be considered as a weaker generalization of the Jung-van der Kulk Theorem for polynomial algebras of rank $n\geq 3$.  In particular, the Nagata automorphism is wild but can be approximated by tame automorphisms \cite{SU25}.

An analogue of the Anick Theorem \cite{Anick} for free metabelian Lie algebras was proven in 1993 by Bryant and Drensky \cite{BD93-2}. They proved that the group of tame automorphisms is dense in the group of all automorphisms of $M_n$ if $n\geq 4$ over fields of characteristic zero. A sharpened version of this result can be found in \cite{KP16}.

Now we prove that the automorphism (\ref{f3}) is almost tame using the ideas of approximation \cite{BD93-2} and the results of the preceding sections. The characteristic-zero assumption in \cite{BD93-2} is not needed in our proof, since we construct the approximation directly.

\begin{lemma}\label{l9}
The automorphism 
\bes
\t=(x_1+[[x_1,x_2],x_1],x_2+[[x_1,x_2],x_2],x_3,\ldots,x_n)
\ees is almost tame. 
\end{lemma}
\Proof We are going to prove that every automorphism of the form  
\bes
(x_1, x_2-\lambda [[x_2,x_3],x_2],x_3-\lambda[[x_2,x_3],x_3],x_4,\ldots,x_n)
\ees
is almost tame. For this it is sufficient to show that 
\bes
P=(x_1+\lambda[x_2,x_3]y_1,x_2,\ldots,x_n) (x_1, x_2-\lambda [[x_2,x_3],x_2],x_3-\lambda[[x_2,x_3],x_3],x_4,\ldots,x_n)\\
=(x_1+\lambda[x_2,x_3]y_1, x_2-\lambda[x_2,x_3]y_2,x_3-\lambda[x_2,x_3]y_3,x_4,\ldots,x_n)
\ees
is almost tame. 

First, following the ideas of \cite{BD93-2,KP16,SU25}, we construct an almost tame automorphism $E$ whose lower-degree part coincides with that of $P$. Consequently, the lower degree of $P^{-1}E$ is at least $4$. 

Consider a Chein automorphism 
\bes
A=(x_1+[x_2,x_3]y_2-\lambda[x_2,x_3]y_1, x_2,\ldots, x_n). 
\ees
Let $\a=(x_1,x_2+\lambda x_1, x_3,\ldots,x_n)$. Then 
\bes
B=A^{\a}=\a A\a^{-1}
=(x_1+[x_2,x_3]y_2+\lambda[x_1,x_3]y_2,x_2-\lambda[x_2,x_3]y_2-\lambda^2[x_1,x_3]y_2, x_3,\ldots,x_n). 
\ees
Direct calculations give that 
\bes
C=(x_1-[x_2,x_3]y_2,x_2,\ldots,x_n)(x_1,x_2+\lambda^2[x_1,x_3]y_2,x_3,\ldots,x_n) B\\
=(x_1+\lambda[x_1,x_3]y_2+g_1y_2^2y_3,x_2-\lambda[x_2,x_3]y_2+g_2y_2^2y_3,x_3,\ldots,x_n), 
\ees
where 
\bes
g_1=-\lambda[x_2,x_3]+\lambda^2[x_1,x_3]-\lambda^2[x_2,x_3]y_2y_3, \ \ g_2=-\lambda^3[x_1,x_3]+\lambda^3[x_2,x_3]y_2y_3, 
\ees
and 
\bes
D=(x_1+[x_2,x_3]y_2,x_2,\ldots,x_n) (x_1,x_2-\lambda^2[x_1,x_3]y_2,x_3,\ldots,x_n) B^{-1}\\
=(x_1-\lambda[x_1,x_3]y_2,x_2+\lambda[x_2,x_3]y_2,x_3,\ldots,x_n)+D'y_2^2y_3,
\ees
where 
\bes
D'=(-\lambda[x_2,x_3]+\lambda^2[x_1,x_3]+\lambda^2[x_2,x_3]y_2y_3, -\lambda^3[x_1,x_3]-\lambda^3[x_2,x_3]y_2y_3,0,\ldots,0).
\ees
Moreover, 
\bes
D^{(23)}=(x_1-\lambda[x_1,x_2]y_3,x_2,x_3-\lambda[x_2,x_3]y_3,x_4,\ldots,x_n)+D'^{(23)}y_2y_3^2
\ees
and 
\bes
E=CD^{(23)}=C(x_1-\lambda[x_1,x_2]y_3,x_2,x_3-\lambda[x_2,x_3]y_3,x_4,\ldots,x_n)+(CD'^{(23)})y_2y_3^2\\
=(x_1+\lambda[x_2,x_3]y_1,x_2-\lambda[x_2,x_3]y_2,x_3-\lambda[x_2,x_3]y_3,x_4,\ldots,x_n)+E'+(CD'^{(23)})y_2y_3^2,
\ees
where
\bes
E'=(h_1,h_2,h_3,0,\ldots,0)
\ees
and
\bes
h_1=g_1y_2^2y_3-\lambda^2[x_2,x_3]y_1y_2y_3-\lambda^2[x_1,x_3]y_2^2y_3-\lambda g_1y_2^3y_3^2+\lambda g_2y_1y_2^2y_3^2,\\
h_2=g_2y_2^2y_3, \ \ h_3=\lambda^2[x_2,x_3]y_2y_3^2-\lambda g_2y_2^2y_3^3. 
\ees

Thus we constructed an almost tame automorphism $E$ whose lower-degree part coincides with that of $P$. 
We have  
\bes
P^{-1}=(x_1-\lambda[x_2,x_3]y_1, x_2+\lambda[x_2,x_3]y_2,x_3+\lambda[x_2,x_3]y_3,x_4,\ldots,x_n)
\ees
and
\bes
P^{-1}E=(x_1,\ldots,x_n)+ Q, 
\ees
where 
\bes
Q=P^{-1}E'+(P^{-1}CD'^{(23)})y_2y_3^2. 
\ees
We have 
\bes
J(P^{-1}E)=I_n+J(Q). 
\ees
Observe that all elements of the matrix $J(Q)$ are divisible by $y_3$, that is, $J(P^{-1}E)\in SL_n(U,y_3U)$   and the elements of the third column of $J(Q)$ are divisible by $y_3^2$. Moreover, all elements of $J(Q)$ have degree great or equal to $4$. By Lemma \ref{l8}, there exists a matrix $R\in S$ such that $R^{-1}J(P^{-1}E)\in SL_n(U,0)=\{I_n\}$, that is, $R^{-1}J(P^{-1}E)=I_n$. Consequently, $J(P^{-1}E)=R\in S$. This means  $P^{-1}E$ is almost tame. Since $E$ is almost tame,  it follows that $P$ and, therefore, $\t$ are almost tame. 
$\Box$

 Recall that the Jacobian matrix of (\ref{f3}) is the matrix  (\ref{f9}). 

\begin{corollary}\label{c3}
The matrix (\ref{f9}) belongs to $S$. 
\end{corollary}

\section{Elimination of a variable }

\hspace*{\parindent} 

In this section, we prove the general form of Lemma \ref{l8} from Section 5 without the restriction $\mathrm{ldeg}(\overline{a}_{in})\geq 3$. In other words, we present a complete proof of Proposition \ref{p6}, which concerns the reduction of the number of variables in matrices from $G$ modulo $S$. We begin by presenting two additional lemmas.

\begin{lemma}\label{l10} There exists a $3$-reducible matrix $C$ in $T$ such that $C_3$ after setting $y_2=y_3=0$ becomes equal to 
\bes
I_3+\lambda y_1E_{31}=I_3+
\begin{bmatrix}
0 & 0 & 0\\
0 & 0 & 0\\
\lambda y_1 & 0 & 0
\end{bmatrix},
\ees
where $\lambda\in K$. 
\end{lemma}
\Proof We use again the language of automorphisms. Set 
\bes
A=(x_1+[x_2,x_3],x_2,\ldots,x_n),\\
 \a=(x_1,x_2+\lambda x_1,x_3,\ldots,x_n),\\
  B=(x_1,x_2-\lambda^2[x_1,x_3],x_3,\ldots,x_n). 
\ees
Then 
\bes
[A,\a]=A\a A^{-1}\a^{-1}=(f_1,f_2,x_3,\ldots,x_n), 
\ees
where 
\bes
f_1=x_1-\lambda[x_1,x_3]-\lambda [[x_2,x_3],x_3], f_2=x_2+\lambda[x_2,x_3]+\lambda^2[x_1,x_3]+\lambda^2[[x_2,x_3],x_3], 
\ees
and 
\bes
B[A,\a]=(g_1,g_2,x_3,\ldots,x_n),
\ees
where 
\bes
g_1=x_1-\lambda[x_1,x_3]-\lambda [[x_2,x_3],x_3]+\lambda^3[x_1,x_3,x_3,x_3],\\
g_2=x_2+\lambda[x_2,x_3]-\lambda^3[x_1,x_3,x_3]+\lambda^2[[x_2,x_3],x_3]-\lambda^4[x_1,x_3,x_3,x_3]. 
\ees
The matrix $C=J(B[A,\a])$ belongs to $T$ by construction, and obviosly, $3$-reducible. Moreover, 
\begin{align*}
C_3=I_3+\begin{bmatrix} -\lambda y_3+\lambda^3y_3^3 & -\lambda^3y_3^2-\lambda^4y_3^3 & 0 \\ 
		-\lambda y_3^2 & \lambda y_3+\lambda^2y_3^2 & 0\\ 
		\lambda y_1+\lambda y_2y_3-\lambda^3y_1y_3^2 & -\lambda y_2 -\lambda^2 y_2y_3+\lambda^3y_1y_3+\lambda^4 y_1y_3^2 & 0 \end{bmatrix}.
\end{align*}

After setting $y_2=y_3=0$ the matrix $C_3$ 
gives $I_3+\lambda y_1E_{31}$.  $\Box$

\begin{lemma}\label{l11} Let 
\bes
I_n+[a_{ij}]\in SL_n(U,y_1U+y_2U)\cap G.
\ees
Then $\mathrm{ldeg}(a_{ij})\geq 2$ for all $1\leq i,j\leq 2$. 
\end{lemma}
\Proof Let 
\bes
A=[a_{ij}]=A_1+\ldots+A_k,
\ees
where all elements of $A_i$ are homogeneous of degree $i$ for all $1\leq i\leq k$ and let $A_1=[b_{ij}]$. Since $a_{ij}\in y_1U+y_2U$ then $b_{ij}\in Ky_1+Ky_2$. Moreover, since $YA=0$ by Lemma \ref{l2}, we get $YA_1=0$. By Lemma \ref{l2} there exist quadratic elements  $f_1,\ldots,f_n\in M_n$ such that 
\bes
A_1=[b_{ij}]=[\d(f_1) \ldots \d(f_n)]. 
\ees
Therefore all components of $\d(f_i)$ belong to $Ky_1+Ky_2$. 

Let 
\bes
f_i=\sum_{k=1}^{n-1} \a_k[x_n,x_k]+f_i', 
\ees
where $f_i'$ does not depend on $x_n$ and $\a_1,\ldots,\a_{n-1}\in K$. If $\a_j\neq 0$ for some $1\leq j\leq n-1$ then
\bes
\frac{\d f_i}{\d x_j}=-\a_jy_n +\frac{\d f_i'}{\d x_j}\notin Ky_1+Ky_2. 
\ees
Consequently, $\a_1=\ldots =\a_{n-1}=0$ and $f_i$ does not depend on $x_n$. The same discussions give that $f_i$ does not depend on $x_3,\ldots,x_n$. Hence $f_i=\lambda_i [x_1,x_2]$ for all $i$.  
Then 
\bes
A_1=\lambda_1F_{121}+\lambda_2F_{122}+\ldots+\lambda_nF_{12n} 
\ees
and the trace $\mathrm{tr}(A_1)$ of $A_1$ is $\lambda_1y_2-\lambda_2 y_1$. 

Recall that the determinat of $I_n+A$ is $1$. The homogeneous part of degree $1$ of the determinant of $I_n+A_1+\ldots +A_k$ is $\mathrm{tr}(A_1)$. Consequently, $\mathrm{tr}(A_1)=0$. The equation $\lambda_1y_2-\lambda_2 y_1=0$ gives that $\lambda_1=\lambda_2=0$. This means that $b_{ij}=0$ and  $\mathrm{ldeg}(a_{ij})\geq 2$  for all $1\leq i,j\leq 2$. $\Box$

\begin{lemma}\label{l12} Let 
\bes
A=I_n+[a_{ij}]\in SL_n(U,\sum_{p=s}^ny_pU)\cap G.
\ees
Then there exists $R\in S$ such that  $R^{-1}A=I_n+[e_{ij}]\in SL_n(U,\sum_{p=s}^ny_pU)\cap G$
 and $e_{1n},e_{2n}, \ldots, e_{n-1 n}\in
\sum_{p=s}^{n-1}y_pU+y_n^2U$.
\end{lemma}
\Proof For any $a\in U$ and any matrix $R$ over $U$ let $\overline{a}$ and $\overline{R}$ be the result of setting $y_s=y_{s+1}=\ldots =y_{n-1}=0$ in $a$ and $R$, respectively.

 By Lemma \ref{l2}, 
\bes
y_1a_{1j}+y_2a_{2j}+\ldots +y_na_{nj}=0, \ 1\leq j\leq n. 
\ees
 Consequently,  
\bes
y_1\overline{a}_{1j}+\ldots y_{s-1}\overline{a}_{(s-1) j}+y_n\overline{a}_{nj}=0, \ 1\leq j\leq n. 
\ees

Since $A=I_n+[a_{ij}]\in SL_n(U,\sum_{p=s}^ny_pU)$ it follows that $a_{ij}\in \sum_{p=s}^ny_pU$ and, consequently, $\overline{a}_{ij}\in y_nK[y_1,\ldots,y_{s-1},y_n]$, that is, $\overline{a}_{ij}=y_n\overline{b}_{ij}$ for some $\overline{b}_{ij}\in K[y_1,\ldots,y_{s-1},y_n]$. Therefore  $\overline{A}=I_n+[\overline{a}_{ij}]=I_n+[y_n\overline{b}_{ij}]$ and the preceding equality gives that 
\bes
y_1\overline{b}_{1j}+\ldots y_{s-1}\overline{b}_{(s-1) j}+y_n\overline{b}_{nj}=0, \ 1\leq j\leq n. 
\ees
If $(\overline{c}_{1n},\ldots,\overline{c}_{(n-1) n})$ is the result of setting $y_n=0$ in $(\overline{b}_{1n},\ldots,\overline{b}_{(n-1) n})$ then from this we get   
\bee\label{f11}
y_1\overline{c}_{1n}+\ldots y_{s-1}\overline{c}_{(s-1) n}=0.
\eee

First we show that there exists a matrix $R\in S\cap SL_n(U,\sum_{p=s}^ny_pU)$ such that 
\bes
\overline{R}=
C=\begin{bmatrix}
I_{n-1} & \begin{matrix}
y_n \overline{c}_{1n}\\
. \\
.\\
.\\
y_n \overline{c}_{(n-1) n}
\end{matrix} \\
\begin{matrix}
0 & \ldots  & 0 
\end{matrix} & 1
\end{bmatrix}.
\ees
 In fact, 
\bes
C=\begin{bmatrix}
I_{n-1} & \begin{matrix}
y_n \overline{c}_{1n}\\
. \\
.\\
.\\
y_n \overline{c}_{(s-1) n}\\
0\\
.\\
.\\
.\\
0
\end{matrix} \\
\begin{matrix}
0 & \ldots  & 0 
\end{matrix} & 1
\end{bmatrix}
\begin{bmatrix}
I_{n-1} & \begin{matrix}
0\\
.\\
.\\
.\\
0\\
y_n \overline{c}_{sn}\\
. \\
.\\
.\\
y_n \overline{c}_{(n-1) n}
\end{matrix} \\
\begin{matrix}
0 & \ldots  & 0 
\end{matrix} & 1
\end{bmatrix}. 
\ees
By (\ref{f11}) and Lemma \ref{l2}, the first factor belongs to $G$. Since it differs from the identity only in column $n$, it represents a Chein automorphism and therefore belongs to $S$. 
 The second factor is equal to the product 
\bes
D=\prod_{k=s}^{n-1} (I_n+y_n\overline{c}_{kn}E_{kn})
\ees
Let 
\bes
\overline{c}_{kn}=y_1f_1+\ldots+y_{s-1}f_{s-1}+\lambda, 
\ees
where $f_i\in K[y_1,\ldots,y_{s-1}]$ and $\lambda\in K$. 
Then 
\bes
y_n\overline{c}_{kn}=y_1y_nf_1+\ldots+y_{s-1}y_nf_{s-1}+y_n\lambda. 
\ees

Let 
\bes
P_k=\begin{bmatrix}
I_{n-1} & \begin{matrix}
-y_ky_n f_1\\
. \\
.\\
.\\
-y_ky_n f_{s-1}\\
0\\
y_1y_nf_1+\ldots+y_{s-1}y_nf_{s-1}\\
0
\end{matrix} \\
\begin{matrix}
0 & \ldots  & 0 
\end{matrix} & 1
\end{bmatrix}\in S
\ees
for any $k$ satisfying $s\leq k\leq n-1$. 

Let $s=n-1$. Then $k=n-1$ and $a_{ij}\in y_{n-1}U+y_nU$. By Lemma \ref{l11}, $\mathrm{ldeg}(a_{ij})\geq 2$ for all $n-1\leq i,j\leq n$. By construction the linear part of $y_n\overline{c}_{kn}$ is obtained from the linear part of $a_{(n-1) n}$. Consequently, $\lambda=0$ and $\overline{P}_k=I_n+y_n\overline{c}_{kn}E_{kn}$. 

Let $s< n-1$, i.e., the set $\{s,\ldots,n\}$ contains at least $3$ elements. 
Then, by Lemma \ref{l10}, there exists a matrix $Q_k\in T$ such that 
\bes
\overline{Q_k}=I_n+\lambda y_nE_{kn}. 
\ees
 Consequently,  $Q_kP_k\in S$ and $\overline{Q_kP_k}=I_n+y_n\overline{c}_{kn}E_{kn}$. 

Therefore in both cases there exists a matrix $\g_i\in S$ such that $\overline{\g_i}=I_n+y_n\overline{c}_{kn}E_{kn}$. Then 
$\overline{\g_s\ldots \g_{n-1}}=D$. Consequently, there exists a matrix $R\in S\cap SL_n(U,\sum_{p=s}^ny_pU)$ such that $\overline{R}=C$. Then $R^{-1}A\in SL_n(U,\sum_{p=s}^ny_pU)\cap G$.
Let $R^{-1}A=I_n+[e_{ij}]$. Recall that, by construction, $\overline{a}_{kn}=y_n\overline{c}_{kn}+y_n^2\overline{d}_{kn}$ for all $1\leq k\leq n-1$, where $\overline{d}_{kn}\in K[y_1,\ldots,y_{s-1},y_n]$. Then 
\begin{align*}
\overline{R^{-1}A}=C^{-1}(I_n+[\overline{a}_{ij}])=\begin{bmatrix}
I_{n-1} & \begin{matrix}
-y_n \overline{c}_{1n}\\
. \\
.\\
.\\
-y_n \overline{c}_{(n-1) n}
\end{matrix} \\
\begin{matrix}
0 & \ldots  & 0 
\end{matrix} & 1
\end{bmatrix}
\left(I_n+\begin{bmatrix}
[\overline{a}_{ij}]_{i,j\leq n-1} & \begin{matrix}
 \overline{a}_{1n}\\
. \\
.\\
.\\
\overline{a}_{(n-1)n}
\end{matrix}\\
\begin{matrix}
	\overline{a}_{n1} & \ldots  & \overline{a}_{n(n-1)} 
\end{matrix} & \overline{a}_{nn} 
\end{bmatrix}\right)\\
=I_n+\begin{bmatrix}
[\overline{a}_{ij}-y_n\overline{c}_{in}\overline{a}_{nj}]_{i,j\leq n-1} & \begin{matrix}
 \overline{a}_{1n}-y_n \overline{c}_{1n}-y_n\overline{c}_{1n}\overline{a}_{nn}\\
. \\
\overline{a}_{kn}-y_n \overline{c}_{kn}-y_n\overline{c}_{kn}\overline{a}_{nn}\\
.\\
 \overline{a}_{(n-1) n}-y_n \overline{c}_{(n-1) n}-y_n\overline{c}_{(n-1)n}\overline{a}_{nn}
\end{matrix}\\
 \begin{matrix}
 	\overline{a}_{n1} & \ldots  & \overline{a}_{n(n-1)} 
 \end{matrix} & \overline{a}_{nn}
\end{bmatrix}
=
I_n+
[\overline{e}_{ij}]
\end{align*}
This implies that every $\overline{e}_{ij}$ is divisible by $y_n$,  since every $\overline{a}_{ij}$ is divisible by $y_n$. Moreover, 
\bes
\overline{e}_{kn}=\overline{a}_{kn}-y_n \overline{c}_{kn}-y_n\overline{c}_{kn}\overline{a}_{nn}=y_n^2\overline{d}_{kn}-y_n\overline{c}_{kn}\overline{a}_{nn}
\ees
is divisible by $y_n^2$. Consequently, $e_{ij}\in \sum_{p=s}^ny_pU$ for all $i,j$ and $e_{kn}\in \sum_{p=s}^{n-1}y_pU+y_n^2U$ for all $1\leq k\leq n-1$.  
The matrix $R^{-1}A$ satisfies the statement of the lemma.  $\Box$

For convenience, we formulate the following complete version of Lemma \ref{l8}.

\begin{lemma}\label{l13}
Let $A=I_n+[a_{ij}]\in SL_n(U,\sum_{p=s}^ny_pU)\cap G$ and  
$a_{1n},a_{2n}, \ldots, a_{(n-1) n}\in
\sum_{p=s}^{n-1}y_pU+y_n^2U$. Then there exists $B\in S$ such that $B^{-1}A\in SL_n(U,\sum_{p=s}^{n-1}y_pU)$. 
\end{lemma}
\Proof  For any $a\in U$, let $\overline{a}$ denote the element obtained from $a$ by setting $y_s=y_{s+1}=\ldots=y_{n-1}=0$. By the assumption of the lemma, $y_n^2$ divides $\overline{a}_{in}$ for all $1\leq i\leq n-1$. Consequently, $\mathrm{ldeg}(\overline{a}_{in})\geq 2$. 

In the proof of Lemma \ref{l8}, the assumption $\mathrm{ldeg}(\overline{a}_{in})\geq 3$ was used only once, namely, to show that 
 the matrix 
\bes
I_n+y_k\overline{d}_{kn}F_{knk}+y_n\overline{d}_{kn}F_{knn}
\ees
belongs to $S$ by Proposition \ref{p2},  since $\mathrm{ldeg}(\overline{d}_{kn})\geq 1$ under this assumption. By Corollary \ref{c3}, this matrix belongs to $S$ even under the weaker assumption $\mathrm{ldeg}(\overline{d}_{kn})\geq 0$. Thus, the proof of Lemma \ref{l8} carries over without the assumption $\mathrm{ldeg}(\overline{a}_{in})\geq 3$ for all $1\leq i\leq n-1$, yielding the desired conclusion. 
$\Box$

\begin{proposition}\label{p6} Let 
\bes
A=I+[a_{ij}]\in SL_n(U,\sum_{p=s}^ny_pU)\cap G.
\ees
Then there exists $B\in S$ such that $B^{-1}A\in SL_n(U,\sum_{p=s}^{n-1}y_pU)$. 
\end{proposition}
\Proof By Lemma \ref{l12}, there exists a matrix $C\in S$ such that $C^{-1}A=I+[c_{ij}]\in SL_n(U,\sum_{p=s}^ny_pU)$ and 
$c_{1n},c_{2n}, \ldots, c_{n-1 n}\in
\sum_{p=s}^{n-1}y_pU+y_n^2U$. By Lemma \ref{l13}, there exists a matrix $D\in S$ such that $D^{-1}C^{-1}A\in SL_n(U,\sum_{p=s}^{n-1}y_pU)$. $\Box$

\section{Main results}

\hspace*{\parindent} 

\begin{theorem}\label{t1}
Every automorphism of the free metabelian Lie algebra $M_n$ of rank $n\geq 4$ over an arbitrary field $K$ is almost tame.  
\end{theorem}
\Proof First we show that $S=G$. Let $A\in G$. Then 
\bes
A\in SL_n(U,\sum_{j=1}^ny_jU). 
\ees
Note that Proposition \ref{p6} allows us to eliminate the last variable $y_n$ from $A$. Since all variables play symmetric roles, the same proposition can be applied to eliminate any variable. Consequently, there exist matrices $B_1,\ldots,B_n\in S$ such that 
\bes
B_s^{-1}\ldots B_1^{-1}A\in SL_n(U,\sum_{j=s+1}^ny_jU)
\ees
for all $s<n$ and 
\bes
B_n^{-1}\ldots B_1^{-1}A\in SL_n(U, 0)=\{I\}.  
\ees
Consequently,  
\bes
B_n^{-1}\ldots B_1^{-1}A=I
\ees
and 
\bes
A=B_1\ldots B_n\in S. 
\ees
This proves that $S=G$. 
Therefore $\mathrm{ATAut}(M_n)\cap \mathrm{IAut}(M_n)=\mathrm{IAut}(M_n)$.  It follows that 
$\mathrm{ATAut}(M_n)=\mathrm{Aut}(M_n)$ since both groups contain all linear automorphisms. $\Box$

\begin{theorem}\label{t2} The automorphism group $\mathrm{Aut}(M_n)$ of the free metabelian Lie algebra $M_n$ of rank $n\geq 4$ over a field $K$ of characteristic different from $3$ is generated by its linear automorphisms together with the quadratic automorphism (\ref{f1}) and the cubic automorphism (\ref{f2}).
\end{theorem}

\begin{corollary}\label{c4} The automorphism group $\mathrm{Aut}(M_n)$ of the free metabelian Lie algebra $M_n$ of rank $n\geq 5$ over a field $K$  is generated by its linear automorphisms together with the quadratic automorphism (\ref{f1}) and the cubic automorphism (\ref{f2}).
\end{corollary}

Theorem \ref{t2} and Corollary \ref{c4} are direct corollaries of Corollary \ref{c1} and Theorem \ref{t1}.

\section*{Acknowledgments}

The author thanks Professor L. Makar-Limanov for useful discussions. 

The project is supported by the grant AP23486782 of the Ministry of Education and Science of the Republic of Kazakhstan.

\end{document}